\documentclass[11pt]{article}
\usepackage[english]{babel}
\usepackage[hmarginratio=3:2]{geometry}
\usepackage[titletoc]{appendix}
\usepackage{authblk}
\usepackage{graphicx}
\usepackage{amsmath}
\usepackage{amsthm}
\usepackage{amssymb}
\usepackage[mathscr]{eucal}
\usepackage{textcomp}
\usepackage[latin1]{inputenc}
\usepackage{latexsym}
\usepackage{hyperref}
\usepackage[dvipsnames]{xcolor}
\usepackage{framed}
\usepackage{tikz-cd}
\usepackage{ifsym}
\usepackage{enumerate}
\usepackage{sectsty}
\usepackage{chngcntr}

\numberwithin{equation}{section}

\numberwithin{subsection}{section}

\allowdisplaybreaks[1]

\newtheorem{theorem}{Theorem}[section]
\newtheorem{lemma}[theorem]{Lemma}
\newtheorem{corollary}[theorem]{Corollary}
\newtheorem{proposition}[theorem]{Proposition}

\newtheorem*{remark}{Remark}
\newtheorem*{definition}{Definition}

\newenvironment{enumerate1}
{\begin{enumerate}[\upshape (1)]}
{\end{enumerate}}

\newcommand\cA{\mathcal{A}} 
\newcommand\cC{\mathcal{C}} 
\newcommand\cE{\mathcal{E}} \newcommand\cF{\mathcal{F}}
\newcommand\cG{\mathcal{G}} \newcommand\cH{\mathcal{H}}
\newcommand\cI{\mathcal{I}} 
 
 \newcommand\cN{\mathcal{N}}
\newcommand\cO{\mathcal{O}}

 \newcommand\cX{\mathcal{X}}

\newcommand\QQ{\mathbb{Q}}

 \newcommand\ZZ{\mathbb{Z}}

\newcommand\arr{\ifinner\to\else\longrightarrow\fi}

\newcommand\arrto{\ifinner\mapsto\else\longmapsto\fi}

\def\displaytimes_#1{\mathrel{\mathop{\times}\limits_{#1}}}

\def\displayotimes_#1{\mathrel{\mathop{\bigotimes}\limits_{#1}}}

\newlength{\ignora}

\newcommand{\mmu}{\boldsymbol{\mu}}

\renewcommand{\epsilon}{\varepsilon}

\renewcommand{\k}[1][*]{\operatorname{K}'_{#1}}

\newcommand{\R}{\mathrm{R}}

\newcommand\nome{testing}
\newcommand\call[1]{\label{#1}\renewcommand\nome{#1}}
\newcommand\itemref[1]{\item\label{\nome;#1}}

\title{On a decomposition theorem in equivariant generalized homology theories for finite group actions}
\author[1]{Francesco Sala}
\affil[1]{\textit{Scuola Normale Superiore, Piazza dei cavalieri 7, 56126, Pisa, Italy}}
\date{}

\begin{document}
\maketitle

\begin{abstract}
In \cite{Vistoli}, A. Vistoli proved a decomposition theorem for the rational equivariant algebraic K-theory of a variety under the action of a finite group $G$. We generalize his result to more general algebraic (co)homology theories having the Mackey property and admitting localization long exact sequences. In general, the pieces are indexed by conjugacy classes of subgroups of $G$.  Our construction is based on some result about a decomposition of the rational Burnside ring of a finite group, which stands behind the classical splitting theorems for equivariant spectra in stable equivariant homotopy theory. 

Applying this result to the case of Borne's modular K-theory we exhibit a case where the splitting is indexed by not necessarily abelian subgroups.
\end{abstract}

\bigskip
\bigskip
\section{The general form of a decomposition theorem for quotient stacks}

Here we discuss a general method for obtaining decomposition theorems with respect to some cohomological theory of a quotient DM stack $\cX=[X/G]$.\\

This discussion should heuristically encompass the well-known decomposition theorems in equivariant algebraic K-theory (by Thomason, Vistoli, Vezzosi, Toen, ...), as well as the decomposition theorems in equivariant homotopy theory (such as Hopkins-Kuhn-Ravenel work about equivariant Morava K-theories and the Minami-Webb formula in equivariant cohomology). \\

In \cite{Vistoli}, A. Vistoli proved a decomposition theorem for the rational equivariant algebraic K-theory of a variety under the action of a finite group $G$. His result states that - over an algebraically closed field -
\[
K(X,G)\simeq \underset{n,\mu_n<G}{\prod} (K(X^\sigma)\otimes \QQ(\zeta_n))^{N_G(\mu_n)}
\]
where the product is indexed by conjugacy classes of cyclic subgroups of $G$. \\

In algebraic topology, there are many similar results: for example the character theory of Hopkins, Kuhn and Ravenel (\cite{HKR}) produces an analogous formula for a special class of cohomology theories, the \emph{complex oriented} ones. In this case, given a theory $E^*$, the decomposition is indexed by the class of abelian subgroups of $G$; the groups $\QQ(\zeta_n)$ must be substituted by the piece of any $E^*(A)$ corresponding to $A$ itself, where $A$ is any abelian group.\\

In this paper we show that in an algebraic context we can prove a result analogous to the classical decomposition of equivariant spectra (see \cite{May}) for a generalized (co)homology theory. The form of our result is analogous to Vistoli's: given an adequate equivariant (co)homology theory $\cH^*(-,G)$ with an action of the rational Burnside ring $B(G)$, we have a decomposition
\[
\cH(X,G)\simeq \prod_{H<G}\cH(X^H,H)^{N_G(H)}_H
\]
where the index $(-)_H$ means that we are localizing with respect to a decomposition of the Burnside ring $B(G)\simeq \prod_{H<G} B(G)_H$.\\

Moreover, we suggest a construction that should settle the question of functoriality of this formula. Namely, following the example of B. Toen (\cite{To}), who gives a map $K([X/G])\rightarrow K(\cI[X/G])$ to the K-theory of the inertia stack, we exhibit a geometric object that in our context plays the same role as the inertia stack for K-theory: it is the \emph{stack of subgroups} (see section 1.2). Then the study of functoriality for the decomposition formula becomes equivalent to the study of functoriality of the stack of subgroups.

\subsection{Conlon's decomposition}

In this section we briefly discuss our main tool, that is Conlon's results about a decomposition of the Burnside algebra and the split K-theory of a finite group. We will follow \cite{conl}, making some slight simplification somewhere.\\

\begin{definition}
Let $G$ be a finite group. The \emph{Burnside algebra} $B(G)$ is the free abelian group over the set of finite $G-$sets, quotiented by relations of the type $[S]+[T]-[S\sqcup T]$. It has a monoidal structure given by $[S]\otimes [T]:=[S\times T]$.
\end{definition}

The Burnside algebra can be identified with the free abelian group generated by sets of the form  $\{ G/\Gamma\}$, where $\Gamma$ ranges over conjugacy classes of subgroups of $G$. Let us call $\langle\Gamma\rangle$ the class of $G/\Gamma$; then the product is given by $\langle\Gamma\rangle\langle\Gamma '\rangle= \bigoplus_{g\in\cC}\langle\Gamma\cap g^{-1}\Gamma g \rangle$, where $\cC:=\Gamma ' \backslash G \slash\Gamma$ is the set of double cosets.\\

It is easy to see the following:

\begin{proposition}
Let $\Lambda$ be a $\ZZ-$algebra such that for any $\Gamma\hookrightarrow G$ the number $\nu^G_\Gamma=|N_G(\Gamma)/\Gamma|$ is invertible in $\Lambda$.

Then the ring map $B(G)_\Lambda\arr\Lambda$ given by sending any $\langle\Gamma\rangle$ to zero for $\Gamma$ a proper subgroup has a canonical section $s_G$.
\end{proposition}

From now on we will always assume that $B(G)$ has coefficients in some $\Lambda$, and we will omit the subscript.

\begin{definition}
Let $_GB(G)$ be the ideal generated by $s_G$; we will refer to it as the \emph{tautological part} of the Burnside ring.
\end{definition}

Conlon proves the following (\cite{conl}, Theorem 2.6)

\begin{proposition}
Let $_HB(G)$ be the extension of $_HB(H)$ along the canonical induction $B(H)\hookrightarrow B(G)$. Then there is a direct sum decomposition $B(G)=\bigoplus_\Gamma\medskip  _\Gamma B(G)$, induced by a complete set of idempotents.

This decomposition is preserved by push-forward and is also invariant by pull-back, in the sense that $Res^G_\Gamma u^G_{\Gamma '}$ is zero unless $\Gamma ' <\Gamma^g$ for some $g$, and in that case $
Res^G_\Gamma u^G_{\Gamma '}$ is a multiple of $\sum_{\Gamma '^g} u^\Gamma_{\Gamma '^g}$, where the $\Gamma '^g$ are all the $\Gamma-$conjugacy classes of $G-$conjugates of $\Gamma '$.

\end{proposition}

The result essentially reduces to the following: for any $\Gamma$, its tautological part is generated by an idempotent $u_\Gamma$, and the induced element $u^G_\Gamma:=Ind^G_\Gamma(u_\Gamma)$ is \emph{almost} idempotent, in the sense that $u^G_\Gamma\cdot u^G_\Gamma=\nu^G_\Gamma\cdot u^G_\Gamma$.

For the sake of completeness we include a revised proof:

\begin{proof}

This is easily shown by induction on the size of $G$ (the base case being trivial). Suppose thus by induction that for any $\Gamma < G$ we have a decomposition of $B(\Gamma)$ induced by a set of orthogonal idempotents $u^\Gamma_{\Gamma '}$ (indexed by conjugacy classes of subgroups $\Gamma '$). Then we set $u^G_\Gamma:=\frac{1}{\nu^G_\Gamma}Ind^G_\Gamma(u^\Gamma_\Gamma)$ for any proper subgroup $\Gamma$, while we set $u^G_G:=1-\underset{\Gamma <G}{\bigoplus} u^G_{\Gamma }$.

Finally we observe that the span of the $u^\Gamma_{\Gamma '}$ contains the induced image of any $B(\Gamma ')$ for all proper $\Gamma '<\Gamma$.\\

We first prove that the $u^G_\Gamma$ form a complete set of idempotents. By the projection and Mackey formulas we have
\[
u^G_\Gamma\cdot u^G_\Gamma=\frac{1}{\nu^G_\Gamma}\cdot Ind^G_\Gamma(\frac{1}{\nu^G_\Gamma}u_\Gamma^\Gamma\cdot Res^G_\Gamma Ind^G_\Gamma u_\Gamma^\Gamma)=\frac{1}{\nu^G_\Gamma}\cdot Ind^G_\Gamma(\frac{1}{\nu^G_\Gamma}u_\Gamma^\Gamma\cdot (\nu^G_\Gamma\cdot u_\Gamma^\Gamma+\sum_{\Gamma '}Ind^\Gamma_{\Gamma '}u'_{\Gamma '}))
\]
for some $u'_{\Gamma '}\in B(\Gamma ')$. Now, the sum $\sum_{\Gamma '}Ind^\Gamma_{\Gamma '}u'_{\Gamma '}$ is in the span of the $u^\Gamma_{\Gamma '}$ and thus has zero product with $u_{\Gamma}^\Gamma$, by the induction hypothesis. This proves that the $u^G_\Gamma$ are idempotents for $\Gamma <G$.

We now want to check that the $u^G_\Gamma$'s are orthogonal. To see it fix two subgroups $\Gamma, \Gamma ' \leq G$. If one of them is $G$ the orthogonality is true by definition if we know the result for proper subgroups; therefore we may suppose without loss of generality that the two subgroups are proper and that $\Gamma '\nless\Gamma$ (meaning that no $G-$conjugate of $\Gamma '$ is a subgroup of $\Gamma$) and proceed as before:
\[
u^G_\Gamma\cdot u^G_{\Gamma '}=\frac{1}{\nu^G_\Gamma}\cdot Ind^G_\Gamma(\frac{1}{\nu^G_{\Gamma '}}u_\Gamma^\Gamma\cdot Res^G_\Gamma Ind^G_{\Gamma '}u_{\Gamma '}^{\Gamma '})=\frac{1}{\nu^G_\Gamma}\cdot Ind^G_\Gamma(\frac{1}{\nu^G_{\Gamma '}}u_\Gamma^\Gamma\cdot (\sum_{\Gamma ''}Ind^\Gamma_{\Gamma ''}u''_{\Gamma ''}))
\]
for some $u''_{\Gamma ''} \in B(\Gamma '')$, where the $\Gamma ''$ represent proper subgroups of $\Gamma $. Again by induction hypothesis the above product is zero. Finally, the combination of the two above results shows that $u^G_G$ is also idempotent, as we wanted.\\

The invariance of this decomposition by push-forward follows trivially from the definitions, while to see the invariance by pull-back we set up another induction. We want to prove that $Res^G_\Gamma u^G_{\Gamma '}$ is a rational multiple of $u^\Gamma_{\Gamma '}$; again, if $\Gamma ' < G$ the above calculation with Mackey formula gives the result: indeed we may write (calling $\Gamma_g:=\Gamma\cap ^g\Gamma '$)
\[
Res^G_\Gamma u^G_{\Gamma '}=Res^G_\Gamma Ind^G_{\Gamma '}u^{\Gamma '}_{\Gamma '}=\underset{g\in \Gamma\backslash G\slash \Gamma '}{\sum} Ind^\Gamma_{\Gamma_g} Res^{\Gamma '}_{\Gamma_g} u^{\Gamma '}_{\Gamma '}.
\]
By the induction hypothesis, the restriction $Res^{\Gamma '}_{\Gamma_g} u^{\Gamma '}_{\Gamma '}$ is zero unless $\Gamma '=\Gamma_g$; in that case the above sum reduces to $\sum_{\Gamma '^g} u^\Gamma_{\Gamma '^g}$.

Finally, the above results imply that any restriction $Res^G_\Gamma u^G_G$ is zero. Indeed the restriction is a ring map and $1\in B(\Gamma)$ is a linear combination of elements that come from the set $\{u^G_{\Gamma '}\}_{\Gamma ' <G}$ by restriction: $1=\sum_{\Gamma '} a_{\Gamma '}Res^G_\Gamma u^G_{\Gamma '}$; this gives readily that $Res^G_\Gamma u^G_G \cdot 1=\sum_{\Gamma '} Res^G_\Gamma (u^G_G\cdot u^G_{\Gamma '})=0$, as we wanted.
\end{proof}

Suppose that we have a cohomological theory $\cH^*(G)$ with an action of the Burnside ring, making $\cH^*(G)$ a graded $B(G)-$module. Then the splitting of $B(G)$ which we described above gives a decomposition $\cH^*(G)=\bigoplus_\Gamma H^*(G)_\Gamma$, where the subscript $-_\Gamma$ indicates the elements with support on $_\Gamma B(G)$.\\

Suppose for example that $\cH^*(G)= K^*(R[G]-\bf{mod}^{split})$ or $\cH^*(G)= K^*(R[G]-\bf{mod})$. In those cases the action of the Burnside algebra is induced by the map $B(G)\arr R[G]-\bf{mod}$ sending $\langle\Gamma\rangle$ to the representation $Ind^G_\Gamma 1$. Conlon (\cite{conl}) shows the following

\begin{theorem}\call{semidirect} Let $R$ be a ring, $\cA=R[G]-\bf{mod}$.
\begin{enumerate1}
\itemref{1} Suppose that  $\cH^*(G)= K^*(\cA^{split})$ or $\cH^*(G)= K^*(\cA)$.

For any $H<G$ the induction and restriction maps
\[
\cH^*(H)_H^{N_G(H)}\xrightarrow{Ind^G_H} \cH^*(G)_H\xrightarrow{Res^G_H} \cH^*(H)_H^{N_G(H)}
\]
are isomorphisms.
\itemref{2} Suppose that $R$ is a complete local ring with residue field of characteristic $p$ and $\cH^*(G)= K^*(\cA^{split})$. Then $\cH^*(G)_H$ is zero unless $H$ is a semidirect product $P\rtimes C$, where $P$ is a $p-$group and $C$ is a cyclic group of order coprime with $p$.
\itemref{3} Suppose that $\cH^*(G)= K^*(\cA)$. Then  $\cH^*(G)_H$ is zero unless $H$ is a cyclic group of order coprime with the characteristic of $R$. 
\end{enumerate1}

\end{theorem}

Note that the first point is not completely trivial, since even though the restriction of $u^G_H$ to $K^*(H)$ is $u_H$ the composition $Res^G_HInd^G_H(a\cdot u_H)$ is not the identity for a general $a$ (indeed it is equal to $u_H\cdot \sum_{\gamma\in N(H)/H}a^\gamma$).

Now, given a suitable cohomology theory $\cH^*(G,X)$, suppose that there is an action of the Burnside algebra $B(G)$ on it. Then, by the prevoius section, there is a corresponding decomposition $\cH^*(G,X)=\bigoplus_H \cH^*(G,X)_H$.

We would like to have some result providing an isomorphism of the $H$-parts $ \cH^*(\cA,X)_H$ with something localized on the fixed-point locus $X^H$. In orded to have a meaningful and functorial definition, we need to explore a variant of the inertia stack construction.

\subsection{The "wild " inertia}

Let us recall given a DM stack $\cX$, there is a classical notion of $n-$inertia
$$\cI_n\cX:=\underline{Hom}^{rep}(B\ZZ\slash n\ZZ,\cX),$$
a stack that classifies representable morphisms from the classifying stack of a cyclic group. Similarly we have a \emph{cyclotomic inertia}, that is defined analogously with $\mu_n$ instead of $\ZZ\slash n\ZZ$.\\

We can generalise this construction. Let us consider, for a finite discrete group $G$, the \emph{universal} classifying stack $B\cG_{[G]}$, that is the classifying stack $BG$ seen over the stack classifying automorphisms of $G$: in other words the stack 
\[
B\cG_{[G]}:=B (G \rtimes (Aut(G)))
\]
with the tautological action of $Aut(G)$ on $G$, which classifies the "twisted" étale group schemes, locally isomorphic to $G$ in the étale topology (see \cite{yas}, Section 3). Over $B\cG_{[G]}$ lives the following object:
.
\begin{definition}
Let us define the stack $I^{[G]} \cX $, so that an object of $I^{[G]} \cX(S)$ is a pair $(\xi,G')\subseteq (\cX\times B\cG_{[G]})(S)$ such that $G'$ is a subgroup of $\underline{Aut}_S(\xi)$. Moreover, let the \emph{wild inertia} be $\cI^{[\cdot]}\cX:= \bigsqcup_G I^{[G]} \cX$.
\end{definition}

For a stack of the form $\cX=[X/G]$, with $G$ an affine group, it is easy to give an explicit description of $\cI^{[\cdot]}\cX$ (see for example \cite{yas}, Lemma 6.16):

\begin{theorem}
Let $\cX=[X/G]$; then we have an isomorphism
\[
I^{[H]} \cX \simeq \underset{H'\in\cC}{\bigsqcup} [X^{H'}/N_G(H')]
\]
where $\cC$ is the set of conjugacy classes of subgroups $H'<G$ such that $H'\simeq H$.
\end{theorem}

\subsection{The general decomposition theorem}

Now, suppose that we have a sufficiently \emph{well-behaved} cohomological theory $\cH^*(G,-)$; in this case we can look for a result of the following form: for each $H<G$ conjugacy class of subgroups the induction 
\[
i_*\colon \cH^*(H,X^H)^{N_G(H)}_H\arr \cH^*(G,X)_H
\]
is an isomorphism.\\

We shall be more precise:

\begin{theorem}
Let $A=Spec(R)$ be a base scheme $\cH^*(-,G)$ be an equivariant (co)homology theory for $A$-schemes with coefficients in $\QQ$ and the following properties:
\begin{enumerate1}
\itemref{1} It admits functorial pull-backs and push-forwards for equivariant maps $(Y,H)\arr (X,G)$.
\itemref{2} It admits an action of the Burnside algebra such that any pull-back or push-forward is a map of $B(G)$-modules.
\itemref{3} It admits localization long exact sequences.
\itemref{4} Depends on $[X/G]$, in the sense that if $X=Y\times_HG$ then $\cH^*(X,G)\simeq \cH^*(Y,H)$.
\itemref{5} When we consider the map $i\colon [X/H]\arr [X/G]$, then the compositions $i^*i_*$ and $i_*i^*$ are consistent with the expectation from representation theory; this means that they satisfy a suitable form of the Mackey formula
\[
i^*i_*(\xi)=\sum_i Ind^H_{H_i^g}\mu^*_g Res^H_{H_i}(\xi)
\]
(where $\mu_g$ indicates the morphism of multiplication by $g\in G$) and of the projection formula
\[
i_*i^*(\xi)=\xi\cdot i_*(1).
\]

\end{enumerate1}

Then by $(1)$ and $(2)$ there is a decomposition $\cH^*(-,G)=\bigoplus_{H}\cH^*(-,G)_H$ indexed by $H<G$ conjugacy classes of subgroups.

We have that, for any $H$, the induction 
\[
i_*\colon \cH^*(H,X^H)^{N_G(H)}_H\arr \cH^*(G,X)_H
\]
is an isomorphism.

\end{theorem}

Moreover, if the action of $B(G)$ is induced, for example, by an action of $\cH^0(G,pt)$, then the only relevant $H$'s are those for which the $H-$parts $\cH^0(G,pt)_H$ are non-zero. This in general leaves out only some restricted class of subgroups.

\begin{proof}
We prove the theorem for a homology theory. The dual case is completely analogous.\\

We begin by observing that any $X$ has an open dense subset $U$ such that $[U/G]$ is a gerbe trivialized by a finite Galois cover. In particular $X$ admits a stratification
\[
U_0\xrightarrow{j_0} U_1\xrightarrow{j_1}\cdots \xrightarrow{j_{n-1}} U_n=X
\]
by open $G-$invariant sets such that for any stratum $X_i:=U_{i+1}\backslash U_i$ the quotient stack $[X_i/G]$ is a gerbe trivialized by a Galois cover with group $\Gamma_i$.\\

By $(3)$ we have fon any $l$ a commuting diagram
\[
\begin{tikzcd}
\cH_{l+1}(U_i^H,H)\ar[r]\ar[d]&\cH_{l}(X_i^H,H)\ar[r,"i_{i,*}"]\ar[d]
&\cH_{l}(U_{i+1}^H,H)\ar[r,"j_i^*"]\ar[d] &\cH_{l}(U_i^H,H)\ar[r]\ar[d]& \cH_{l-1}(X_i^H,H)\ar[d]\\
\cH_{l+1}(U_i,G)\ar[r]&\cH_{l}(X_i,G)\ar[r,"i_{i,*}"]
&\cH_{l}(U_{i+1},G)\ar[r,"j_i^*"] &\cH_{l}(U_i,G)\ar[r]& \cH_{l-1}(X_i,G)
\end{tikzcd}
\]
From the five lemma we immediately infer that if we know the thesis for $U_i$ and for any $X$ such that $[X/G]$ is a gerbe Galois-locally trivial then we also have it for $U_{i+1}$.\\

All we have to do, then, is to show the thesis in the case when $[X/G]$ is a gerbe trivial on a Galois cover of the moduli space $M$. 

Let $\Gamma$ be the Galois group of the trivializing cover $N\arr M$ and $H$ be the stabilizer. Let $p$ be the projection $[N/H]\arr [X/G]$. Then by $(4)$ and $(5)$ we immediately deduce that the pull-back $p^*$ induces an isomorphism 
\[
\cH^*(X,G)\simeq \cH^*(N,H)^\Gamma
\]  
Indeed, the compositions $p^*p_*$ and $p_*p^*$ are easily seen to be equal to the multiplication by $|\Gamma|$ map.\\

We are left with proving the statement for a trivial gerbe $[X/G]$, where $G$ acts trivially. For that we just repeat the proof of Proposition~\ref{semidirect}. This is actually quite immediate from $(5)$: first of all, by the projection formula we have $Ind^G_HRes^G_H(\xi)=\xi\cdot Ind^G_H(1)$ and $Ind^G_H(1)$ is invertible in $B(G)_H$, making the composition an isomorphism; moreover, by the Mackey formula, we have $Res^G_HInd^G_H(\xi)=\sum_i Ind^H_{H_i^g} Res^H_{H_i}(\xi)$; but when $\xi$ is in the $H-$part its restriction to any proper subgroup of $H$ is zero, so that we can rewrite $Res^G_HInd^G_H(\xi)=\sum_{n\in N_G(H)/H} n\cdot\xi=|N_G(H)/H|\cdot\xi$ (where the action of $N_G(H)/H$ is induced by the conjugation action on $H$) because $\xi$ is $N_G(H)-$invariant. We conclude that in this case the induction and restriction maps
\[
\cH_*(H)_H^{N_G(H)}\xrightarrow{Ind^G_H} \cH_*(G)_H\xrightarrow{Res^G_H} \cH_*(H)_H^{N_G(H)}
\]
are isomorphisms, and the thesis follows.
\end{proof}

\subsection{The decomposition of algebraic K-theory of DM stacks}

Suppose for instance that $\cH_*(X,G)=\k([X/G])$. K-theory is well-known to satisfy all the properties of the Theorem, so we conclude that the result holds for the K-theory of quotient Deligne-Mumford stacks.\\

If the action of $G$ on $X$ is tame and the base ring contains all roots of unity we can make a comparison with the more classical decompositions (see e.g. \cite{Vistoli}). In that case we have a decomposition of the representation ring of free $G$-representations $\R G=\underset{\mmu_r\in G}{\bigoplus} \QQ(\zeta_r)^{N_G(\mmu_r)}$ indexed by conjugacy classes of cyclic subgroups of $G$; for any $H\simeq \mmu_r$ let us call $\R G '_H$ the $H-$piece in this classical decomposition. 

\begin{theorem}
 The classical decomposition coincides with the one induced by the Burnside algebra.
\end{theorem}

\begin{proof}
This is immediate to see that when $G$ is itself cyclic, by induction on the dimension of $G$. More precisely, we know by Theorem ~\ref{semidirect} that Conlon's decomposition of $\R G$ is also indexed by conjugacy classes of cyclic subgroups, and that for any $H$ we have an isomorphism $\R G_H\simeq \R H_H^{N_G(H)}$ induced by the induction $\R H\arr \R G$; thus it sufficies to see that $\R G_G$ coincides with $\R G'_H$. But we may suppose that we know it for all the smaller groups (by induction hypothesis); if we call  $v^G_H $ the idempotent generating $\R G_H '$, then $\R G_H '$ is generated by the only rational multiple of $Ind^G_H v^H_H $ which is itself idempotent, so it is uniquely determined. In particular this generator coincides with $u^G_H$, by inductive hypothesis: this implies that $u^G_G=1-\sum_H u^G_H=1-\sum_H v^H_H =v^G_G $, as we wanted.

Finally, the covariance by push-forward immediately settles the general case.
\end{proof} 

In particular the classical decomposition $\k([X/G])=\bigoplus_H \k([X/G])_H$  coincides with the one induced by the Burnside ring, since they are both obtained localizing $\k([X/G])$ with respect to this common $\R G$-decomposition.\\

Now, using the self-intersection formula, we can give inverse maps
\[
\k(X,G)_H\arr \k(X^H,H)^{N(H)}_H.
\] 

First of all, by Proposition~\ref{semidirect} the only relevant $H$'s are the cyclic groups. Given one such $H$, let now $\cN$ be the conormal sheaf associated to the closed immersion $i\colon X^H\hookrightarrow X$.

Suppose that $i$ is a regular immersion. By the self-intersection formula, we know that the composition 
\[
\k(X^H,H)^{N(H)}_H \xrightarrow{i_*} \k(X,H)^{N(H)}_H
\xrightarrow{i^*} \k(X^H,H)^{N(H)}_H
\]
is equal to the multiplication by the element $\lambda_{-1}(\cN)$. But we have
\begin{lemma}
The element $\lambda_{-1}(\cN)$ is invertible in $\k(X^H,H)_H$.
\end{lemma}

\begin{proof}
Let $Y:=X^H$. If $Y$ does not have points in characteristic not coprime with $|H|$ then the result is known from the classical theory (see \cite{To} or \cite{Vistoli}). 

Otherwise, let $Y'\hookrightarrow Y$ be the closed subscheme of points whose residue characteristic is not coprime with $|H|$. By Proposition~\ref{semidirect} we have that $\k(Y',H)_H=0$; by the localization exact sequence, naming $j\colon U=(Y\backslash Y')\hookrightarrow Y$, we have that $j^*\colon \k(Y,H)_H\arr \k(U,H)_H$ is an isomorphism. But $j^*$ clearly preserves the multiplication by $\lambda_{-1}(\cN)$, so the thesis follows.
\end{proof}
Let us consider the composition
\[
\k(X^H,H)^{N(H)}_H \xrightarrow{i_*} \k(X,H)^{N(H)}_H\xrightarrow{Ind^G_H}\k(X,G)_H\xrightarrow{Res^G_H}\k(X,H)^{N(H)}_H
\xrightarrow{i^*} \k(X^H,H)^{N(H)}_H.
\]
By Mackey formula the central composition $\k(X,H)^{N(H)}_H\xrightarrow{Ind^G_H}\k(X,G)_H\xrightarrow{Res^G_H}\k(X,H)^{N(H)}_H$ sends $\xi$ to $\sum_i Ind^H_{H_i^g} Res^H_{H_i}(\xi)$; but when $\xi$ is in the $H-$part its restriction to any proper subgroup of $H$ is zero, so that we can rewrite $Res^G_HInd^G_H(\xi)=\sum_{n\in N_G(H)/H} n\cdot \xi=|N(H)/H|\cdot\xi $. We conclude that the composition is equal to the multiplication by the element $|N(H)/H|\cdot \lambda_{-1}(\cN)$. In particular the modified restriction map
\[
\frac{i^*}{|N(H)/H|\cdot \lambda_{-1}(\cN)}\colon \k(X,G)_H\arr \k(X^H,H)^{N(H)}_H
\]
is an isomorphism inverse to $i_*$.\\

Summarizing, we have the following
\begin{theorem}
Let $X$ be an $R-$scheme with an action of a finite constant group $G$. Then we have a decomposition $\k(X,G)=\underset{H\hspace{1mm} \emph{cyclic}}{\bigoplus}\k(X,G)_H$ and the push-forward induces isomorphisms
\[
i_*\colon \k(X^H,H)^{N(H)}_H\arr \k(X,G)_H.
\]
If the immersion $X^H\hookrightarrow X$ is regular there is an inverse
\[
\frac{i^*}{|N(H)/H|\cdot \lambda_{-1}(\cN)}\colon \k(X,G)_H\arr \k(X^H,H)^{N(H)}_H.
\]
\end{theorem}

\section{The modular K-theory of a quotient}

We will now discuss a closely related example in equivariant algebraic geometry, that of modular K-theory.

In this section we briefly recall the theory of modular K-theory, along the lines of the original article \cite{bor} by Niels Borne. Here we will suppose that $X$ is defined over a field $k$.

\subsection{Sheaves of modules over a noncommutative algebra}
We first recall the definitions: suppose that $\cA$ is a subcategory of the category of finite dimensional G-representations $\cA_{tot}$.

Let $ \bf{Mod}( \cA)$ be the category $[\cA^{op}, \bf{Ab}]$ of contravariant functors from $\cA$ to the category of abelian groups (with natural transformations as morphisms). A (right) $\cA$-module
is by definition an object in $ \bf{Mod}( \cA)$. $\cA$-modules obviously form an abelian category.

The Yoneda embedding $\cA \arr  \bf{Mod}( \cA)$ sends an object $V$ to the contravariant representable functor $\underline{V} = \cA(\cdot,V)$. An $\cA$-module is said to be of finite type if it is a
quotient of a finite direct sum of representable functors, and we will call $ \bf{mod}( \cA)$ the full
subcategory of $ \bf{Mod}( \cA)$ consisting of objects of finite type.

\begin{definition}
The projective completion $Q\cA$ of $\cA$ is the full subcategory of $ \bf{Mod}( \cA)$ whose objects are the projective $\cA$-modules of finite type ( direct summands of finite direct sums of representable functors).
\end{definition}

\begin{definition}
We say that $\cA$ admits a finite set of additive generators if it is Morita equivalent to the full subcategory generated by a
finite set of its objects.
\end{definition}

Now, let $\cA$ be a projectively complete full subcategory of $\cA_{tot}$ admitting a finite set of additive generators. Starting from $\cA$ we can form two distinct exact categories: the first is $\cA^{split}$, that is $\cA$ with split exact sequences of $G$-modules; the second is $\bf{mod}( \cA)$ with its natural exact structure induced by the abelian category structure. 

The Yoneda embedding induces a morphism of exact categories $\cA^{split}\arr \bf{mod}( \cA)$, and thus also a map in K-theory $\k(\cA^{split})\arr \k(\bf{mod}( \cA))$. We can formulate (see \cite{bor}, section $2.2$) the following

\begin{theorem}
\begin{enumerate1}

\itemref{1} The two groups $\k(\cA^{split})$ and $\k(\bf{mod}( \cA))$ are abelian free of the same rank.
 
\itemref{2} If $G$ is of finite representation type and $\cA=\cA_{tot}$ then the Yoneda embedding is an isomorphism. 

\end{enumerate1}
\end{theorem}

In particular we see that the K-theory of $\bf{mod}( \cA)$ is a new invariant that we can associate to the category $\cA_{tot}=k[G]-\bf{mod}$ and it is in general finer than the classical algebraic K-theory.\\

We can generalize this discussion to the case of a $G-$scheme $Y$, so that what we have done so far corresponds to $Y=Spec(k)$.

\begin{definition}
Let X be a scheme.

\begin{enumerate1}
\itemref{1}  A ring (with several objects) on Y is, by definition, a category $\cA$ enriched on $ \bf{Qcoh}(Y)$. The pair $(Y, \cA)$ is called a ringed scheme.
\itemref{2}   A morphism of ringed schemes $(Y',\cA')\arr (Y,\cA)$ is a couple $(f,f^{\sharp})$ where $f\colon Y'\arr Y$ is a scheme morphism and $ f^{\sharp}\colon \cA\arr f_*\cA'$ is a morphism of  $\bf{Qcoh}(Y)$-categories.

\itemref{3} A quasicoherent sheaf on $(Y,\cA)$ is by definition an enriched functor from $\cA^{op}$ to $\bf{Qcoh}(Y)$.

\itemref{4} A morphism of $\cA$-sheaves is an enriched natural transformation.

\end{enumerate1}
We denote by $\bf{Qcoh}(Y,\cA)=[\cA^{op},\bf{Qcoh}(Y)]$ the corresponding enriched category, and by $Qcoh(Y,\cA)$ the underlying category.
\end{definition}

This definition of ringed scheme allows for a notion of functoriality, which is essentially induced by the functoriality of  $Qcoh(Y)$.

Consider first the case of a morphism $(1, i):
(Y,\cA')\arr (Y,\cA)$. In this case it is easy to see that left and right Kan extensions along $\cA'\arr \cA$ exist and give left and right adjoints to the restriction functor  $\bf{Qcoh}(Y,\cA) \arr \bf{Qcoh}(Y,\cA')$. The left adjoint
will be denoted by $\otimes_ {\cA'} \cA$.

We then extend this definition to the general case:

\begin{definition}
Let us be given a morphism of ringed schemes $(Y',\cA')\arr (Y,\cA)$.
\begin{enumerate1}
\itemref{1} We define $f_\triangle: f_*\bf{Qcoh}(Y',\cA') \arr \bf{Qcoh}(Y,\cA)$ as the $\bf{Qcoh}(Y)$ functor making
the following diagram commute:
\[
\begin{tikzcd}
   f_*\lbrack\cA'^{op},\bf{Qcoh}( Y')\rbrack \ar[r,"f_*"]
   \ar[dr,"f_{\triangle}",swap] &
   \lbrack f_*\cA'^{op},f_*\bf{Qcoh}( Y')\rbrack
     \ar[d," f^{\sharp} "]\\
    & \lbrack\cA^{op}, \bf{Qcoh}( Y)\rbrack
\end{tikzcd}
\]

\itemref{2} We define $f^\triangle \colon \bf{Qcoh}(Y,\cA) \arr f_*\bf{Qcoh}(Y',\cA')$ as the $\bf{Qcoh}(Y)$ -functor
making the following diagram commute:
\[
\begin{tikzcd}
\lbrack\cA '^{op},\bf{Qcoh}(Y')\rbrack & \\
 \lbrack f^*\cA^{op}, \bf{Qcoh}(Y')\rbrack
 \ar[u,"\otimes_{f^*\cA}\cA '" ] &         \\
   \lbrack f^*\cA^{op}, f^*\bf{Qcoh}( Y)\rbrack
   \ar[u,"f^*"]
   &  f^*\lbrack\cA^{op}, \bf{Qcoh}( Y)\rbrack \ar[l,"f^*"]
   \ar[uul,"adj(f^{\triangle})",swap]
\end{tikzcd}
\]
\end{enumerate1}
\end{definition}

These functors satisfy the natural adjunction requirement:

\begin{proposition}
The couple $(f^{\triangle},f_\triangle)$ is part of a
$\bf{Qcoh}( Y)$-adjunction between $\bf{Qcoh}(Y,\cA)$ and $f_*\bf{Qcoh}(Y',\cA')$.
\end{proposition}

For each object $V$ of $\cA$, we denote by $\langle V \rangle$ the
full subcategory of $\cA$ containing only the object $V$.

\begin{proposition}
\begin{enumerate1}
\itemref{1} The category ${Qcoh}(Y,\cA)$ is an abelian category.

\itemref{2}  A sequence $\cF' \arr \cF \arr \cF''$ of $\cA$-sheaves on $Y$ is
exact if and only if for each object $V$ of $\cA$ the sequence
$\cF'(V) \arr \cF(V) \arr \cF''(V)$ is exact in ${Qcoh}
(Y, \langle V \rangle)$.
\end{enumerate1}
\end{proposition}

Since localization commutes with evaluation at each $V$, we immediately get

\begin{corollary}
A sequence of $\cA$-sheaves
$\cF' \arr \cF \arr \cF''$ is exact if and only if for each $p\in Y$ the sequence of stalks
$\cF'_p \arr \cF_p \arr \cF''_p$
is exact.
\end{corollary}

\subsection{Sheaves over the Auslander algebra}

Now we are ready to return to our $G-$scheme $X$. We have an equivalence between the category of equivariant $G-$sheaves on $X$, ${QCoh}(X,G)$, and the category of sheaves on the quotient stack $\cX$. We will denote by $\pi\colon \cX\arr Y$ the map from $\cX$ to the moduli space $Y$ (that is, the quotient scheme).
In particular the category $QCoh(X,G)$ is thus naturally enriched over $QCoh(Y)$ by $\bf{QCoh}(\cF,\cF'):= \pi_*QCoh_\cX(\cF,\cF')$.

\begin{definition}

Let $(X,G)$ be as above, $s_X:X \arr Spec( k)$ be the structure morphism, and $\cA$ a category equipped with a functor $F\colon \cA\arr k[G]-\bf {mod}$.
\begin{enumerate1}
\itemref{1} The \emph{Auslander algebra over $Y$} associated to $(\cA,F)$ is the ring $\cA_X$
over $Y$ equal with the same objects of $\cA$ and morphisms
of the form $\bf{Qcoh}(G,X)(s_X^* F(V),s_X^* F(V'))$, for all objects $V,V'$ of $\cA$.
\itemref{2} The \emph{constant Auslander algebra over $Y$} is the ring
$\cA_X^c=s_Y^*\cA$ (that is, any morphism $Hom_{\cA_X^c}(V,V')$ is given by $QCoh_Y(s_Y^*F(V),s_Y^*F(V'))$).
\itemref{3} The \emph{ "free" Auslander algebra over $Y$} is the ring $\cA_X^f$ over $Y$ with the same objects of $\cA$ and such that any morphism $Hom_{\cA_X^f}(V,V')$ is given by $QCoh_Y(\pi_*s_X^*F(V),\pi_*s_X^*F(V'))$).
\end{enumerate1}
\end{definition}

Clearly $\cA_X\simeq \cA_X^c$ when the $G$-action is trivial and $\cA_X\simeq\cA_X^f$ when the $G-$action is free.\\

The Auslander algebra is functorial:

\begin{proposition}
\begin{enumerate1}
\itemref{1} Any map $f : X' \arr X $ induces a morphism of ringed schemes $(\tilde{f},
  \tilde{f}^\sharp)\colon(Y', \cA_{X'}) \arr (Y, \cA_{X})$.

\itemref{2}  $adj(\tilde {f}^\sharp) : \tilde{f}^*  \cA_{X} \arr \cA_{X'}$
  is an isomorphism if $X' = X \times_Y Y'$.
\end{enumerate1}

\begin{proof}
\begin{enumerate1}
\itemref{1} To construct $\tilde {f}^\sharp$, we start from the isomorphism
  $f^*\bf{Qcoh}(X)({s_X}^*V, {s_X}^*W) \arr \bf{Qcoh}( X')({s_{X'}}^*V, {s_{X'}}^*W)$. This is in fact a $G$-isomorphism and give by adjunction
  a $G$-morphism
  $\bf{Qcoh}( X)({s_X}^*V, {s_X}^*W) \arr f_*\bf{Qcoh}( X')({s_{X'}}^*V,
  {s_{X'}}^*W)$. Applying $\pi_*$ and using the fact that
  $\pi_* f_* = \tilde{f}_*{\pi '}  $, we get the map $\cA_X
  ({s_X}^*V, {s_X}^*W) \arr \tilde{f}_* \cA_{X'} ({s_{X'}}^*V, {s_{X'}}^*W)$
  that we needed.
\itemref{2} By base change, the canonical $2$-arrow $\tilde{f}^* \pi_*
  \Longrightarrow {\pi '}_* f^*$ is an isomorphism, and the result
  follows.
\end{enumerate1}
\end{proof}
\end{proposition}

\begin{definition}
We define $\bf{Qcoh}(\cA,X)$ as the category of quasi-coherent sheaves for the ringed scheme $(Y,\cA_X)$. For any $V\in\cA$ call $\cA_V$ the corresponding representable sheaf.
\end{definition}

There is an obvious map from classical equivariant sheaves to $\cA-$sheaves:
\begin{eqnarray*}
   U_{\cA}\colon\bf{Qcoh}( G,X) \arr &  \bf{Qcoh}(\cA,X) & \\
          \cF        \arr &  \underline{\cF} = (V \arr & \bf{Qcoh}(G, X)(s_X^*V,\cF))
\end{eqnarray*}

We now briefly discuss the functorial properties of this definition.\\

First of all, given an equivariant morphism $f\colon X'\arr X$, it produces an adjoint pair $(f^\triangle, f_\triangle)\colon \tilde{f}_*\bf{Qcoh}(\cA,X')\arr \bf{Qcoh}(\cA,X)$, which we rename as $(f^\cA,f_\cA)$.\\

Moreover, as shown by N. Borne, there are natural induction and restriction functors from and to subcategories:

\begin{proposition}\call{restriction}
Let $\cA'$ be a full subcategory of $\cA$.
\begin{enumerate1}
\itemref{1} The restriction functor $R:\bf{Qcoh}(\cA,X)\arr \bf{Qcoh}(\cA',X) $ admits as a right adjoint the functor $K$
defined on objects by :
\begin{eqnarray*}
  K\colon \bf{Qcoh}(\cA',X) \arr  &  \bf{Qcoh}(\cA,X) & \\
          \cF   \arr &  (V \arr & \bf{Qcoh}(\cA',X)
(\cA_X(\cdot,V)_{|\cA'_X },\cF))
\end{eqnarray*}
Moreover $RK \simeq 1$ (equivalently, $K$ is fully faithful).
\itemref{2} If $\cA$ is projectively complete,
and every object of $\cA$ is a direct summand of a
finite direct sum of objects of $\cA'$, this adjunction is an
equivalence $\bf{Qcoh}(\cA,X) \simeq  \bf{Qcoh}(\cA',X)$.
\end{enumerate1}
\end{proposition}

In particular, when $\cA$ is projectively complete we can assume that $\cA$ is generated by a single object.\\

Finally, with our definitions, there is an obvious notion of change of group: indeed for any morphism $\alpha\colon H\arr G$ we can form the "restricted" functor $\alpha^*\colon \cA\arr k[G]-\bf {mod}\xrightarrow{\alpha^*} k[H]-\bf {mod} $. Naming $Z:=X/H$, there are a morphism $\tilde\alpha\colon Z\arr Y$ and a natural map $\alpha^\sharp\colon \cA_X\arr \tilde\alpha_*(\alpha^*\cA)_X$, giving a morphism of ringed schemes 
\[
(\tilde{\alpha},{\tilde{\alpha}}^{\sharp}): (Z,(\alpha^*\cA)_X)\arr (Y,
\cA_X)
\]

Thus we have restriction functor
$\alpha^{\cA} \colon \bf{Qcoh}(\cA, X) \arr
 \tilde{\alpha}_* \bf{Qcoh}(\alpha^*\cA, X_{|H})$
and an induction functor
$\alpha_{\cA} \colon \tilde{\alpha}_*\bf{Qcoh}(\alpha^*\cA, X_{|H})
\arr \bf{Qcoh}(\cA, X)$. Again, this gives an adjoint pair $(\alpha^\cA,\alpha_\cA)$.\\

Let us make some comment on these induction and restriction functors. Following \cite{bor}, Lemma 7.11, we see that when $G$ is of finite representation type, $X=Spec(k)$ and $\cA= k[G] - \bf{mod}$ then (identifying $k[G]-\bf{mod}$ with a fully faithful subcetegory of $\cA-\bf{mod}$ via the Yoneda embedding) the induction and restriction functors coincide with the usual induction and restriction for group representations.

More precisely, for any subgroup $H\hookrightarrow G$ we have an isomorphism $k[H]-\bf{mod}\simeq \alpha^*\cA$ (this follows immediately by the prevoius proposition and Mackey formula); there is a commuting diagram
\[
\begin{tikzcd}
k\lbrack H\rbrack-\bf{mod}\ar[r,"Ind_H^G"]\ar[d] & k\lbrack G\rbrack-\bf{mod}\ar[d]\\
\alpha^*\cA-\bf{mod} \ar[r,"\alpha_\cA"] & \cA-\bf{mod}
\end{tikzcd}
\]
and a similar one for restriction.

In general, from the induction-restriction (co-)adjunction it is immediate to see that $\alpha_\cA$ (resp. $\alpha^\cA$) sends a representable sheaf $\cA_V$ to the functor $Hom_G(-,Ind^G_H(V))$ (resp. $Hom_H(-,Res^G_H(V))$ \emph{if the latter is representable}.\\

Consider now the case where $\cA$ is generated by a single element $I$; in this case $\cA_X$ can be identified with an actual (non-commutative) algebra $\cA_X=(s_X^* End_k(I))^G$. Let us call $A:=s_X^* End_k(I)$; then for any $\cF$ in $\bf{Qcoh}(\cA,X)$ the induced functor $\alpha_\cA \cF$ is $\alpha_\cA\cF(I)= \tilde\alpha_*\cF(I)$, seen as a $A^H$-module. Conversely the restricted functor is $\alpha^\cA\cF(I)=\tilde\alpha^*F(I)\otimes_{\tilde\alpha^*A^G} A^H$.\\

We end this section with a description of monoidal structures. Suppose that $\cA$ is a full submonoidal category of $k[G] \bf{mod}$, and $X$ is a $G$-scheme.
Then there is a canonical structure of monoidal closed symmetric
category on $\bf{Qcoh}(\cA,X)$, whose unit object is
\[
 \underline{\cO_X} \colon V \arr \pi_*^G(s_X^*V^{\vee})
\]
and whose tensor product is given by:

\[
\cF \otimes \cG (V) = \int^W  \cF(W) \otimes_{\cO_Y} \cG(V\otimes_k W^{\vee})
\]
Tensor product with a representable sheaf can be easily calculated:
\[
 (\cF \otimes \cA_V)(W) \simeq \cF ( W\otimes_k V^{\vee}).
\]
With this it is immediate to see that pull-backs and push-forwards preserve the product with any representable sheaf.

\subsection{Modular K-theory}
As for classical $G-$sheaves, there is a notion of coherence for $\cA_X-$sheaves:
\begin{definition}
A quasicoherent
$\cA$-sheaf $\cF$ on $X$ is said \emph{coherent} if for each
$G$-invariant
open affine $U = Spec(R)$ of $X$, the restriction $\cF_{|U}$ is
of finite type in $\bf{Qcoh}(\cA, U) \simeq [\cA_X(U)^{op}, R^G-\bf{Mod}]$ (i.e. if, seen in $[\cA_X(U)^{op}, R^G-\bf{Mod}]$, it is a quotient of a finite sum of representable objects).

We denote by $Coh(\cA, X)$ the whole subcategory of $Qcoh(\cA,
X)$ whose objects are the coherent $\cA$-sheaves.
\end{definition}

\begin{proposition}
Let $\cA'$ be a full subcategory of $\cA$.
Suppose that $\cA$ is projectively complete,
and that every object of $\cA$ is a direct summand of a
finite direct sum of objects of $\cA'$. Then restriction along
$\cA' \arr \cA$ induces an
equivalence $Coh(\cA,X) \simeq  Coh(\cA',X)$.
\end{proposition}

\begin{proposition}
Suppose $\cA$ admits a
 finite set of additive generators.
Then $Coh(\cA, X)$ is an abelian category.
\end{proposition}

With this notion at our service it is immediate to define a K-theory functor:

\begin{definition}
Let $X$ be a $G$-scheme over $k$ and $\cA$ a full
subcategory of $k[G]-\bf{mod}$  admitting
a finite set of additive generators.
We denote by $K_i(\cA, X)$ the Quillen $i$-th group of the abelian
category $Coh(\cA, X)$. 
\end{definition}

We will refer to this functor as \emph{modular K-theory}.

It admits well-defined notions of pull-back and push-forward:

\begin{proposition}
Let $f: X'\arr X$ be a morphism of $G$-schemes over $k$ such that the
morphism $\tilde{f} \colon Y' \arr Y$ between quotient schemes is flat, and
$X'=X\times_Y Y'$.

Then the functor $f^{\cA}\colon Coh(\cA, X) \arr Coh(\cA, X')$ is
exact, hence induces a map in $K$-theory.
\end{proposition}
\begin{proposition} Let $f: X'\arr X$ be a morphism of $G$-schemes over $k$ and $\cA$ a full
subcategory of $k[G]- \bf{mod}$, admitting a
finite set of additive generators. Suppose that the morphism $\tilde f
\colon Y' \arr Y$ between quotient schemes is proper, then:
\begin{enumerate1}

\itemref{1} For each coherent $\cA$-sheaf $\cF$ on $X'$, and each nonnegative
integer $i$, the $\cA$- sheaf $R^if_{\cA}\cF$ is coherent.

\itemref{2} There exists an integer $n$, such for any integer $i>n$, and any
coherent $\cA$-sheaf $\cF$ on $X'$, we have $R^if_{\cA}\cF=0$.

\itemref{3} If $\tilde f$ is finite or $Y'$ admits an ample line bundle then there is a well-defined push-forward $f_\cA\colon K_i(\cA,X')\arr K_i(\cA,X)$.
\end{enumerate1}
\end{proposition}

It is immediate to see that these definitions are functorial.\\

We conclude this section citing a very important property of modular K-theory, that is it satisfies a form of localization (\cite{bor}, Theorem 6.7):

\begin{theorem}\call{localization}
Let $i\colon X'\arr X$ be a morphism of $G$-schemes over $k$, and $\cA$ a full
subcategory of $k[G] - \bf{mod}$, admitting a
finite set of additive generators.

Suppose that the morphism $\tilde i \colon Y' \arr Y$
between quotient schemes is a closed immersion, and that $i^\sharp \colon
\cA_X \arr \tilde i_*\cA_{X'}$ is an epimorphism \bigskip (*).

Denote by $U$ the pullback by $\pi:X \arr Y$ of the complement of $Y'$
in $Y$, and by $j:U\arr X$ the canonical inclusion.

Then there is a long exact sequence:
\[
\dots \arr K_{i+1}(\cA,U)\arr K_i(\cA,X')\arr K_i(\cA,X)\arr K_i(\cA,U)\arr K_{i-1}(\cA,X')\arr\dots
\]
\end{theorem} 

The application of this theorem is limited by the condition of $i^\sharp \colon
\cA_X \arr \tilde i_*\cA_{X'}$ being an epimorphism. However in some case we can ensure that; for instance, when the substack $[X'/G]\hookrightarrow [X/G]$ is a gerbe (\cite{bor}, Proposition 6.13):

\begin{proposition}\call{gerbe}
Let $i\colon X'\arr X$ be a morphism of $G$-schemes over $k$, such that there is a normal subgroup $H$ of $G$ that acts trivially on $X'$, and $G/H$ acts freely on $X'$; suppose moreover that the morphism $\tilde i \colon Y' \arr Y$ between quotient schemes is a closed immersion.

Then the canonical morphism $i^\sharp \colon \cA_X \arr \tilde i_*\cA_{X'}$ is an epimorphism.
\end{proposition}

\section{A localization theorem for modular K-theory}

Suppose now that $X$ is a quasi-projective variety over an algebraically closed field, equipped with an action of a finite group $G$.

Let $\cA$ be a finitely generated projectively complete full subcategory of $k[G]-\bf{mod}$ containing the trivial representation and equipped with an action of the Burnside algebra (that is, for any $V$ in $\cA$ and any $\Gamma < G$ the tensor product $V\otimes Ind^G_\Gamma 1$ is also in $\cA$).\\

Then we also have an action of the Burnside algebra on $Coh(\cA,X)$, induced by the internal tensor product with representable sheaves. It is immediate to see that this action passes to K-theory, making $K^*(\cA,X)$ a $B(G)$-module. In particular we have a decomposition $K^*(\cA,X)=\bigoplus_{H<G} K^*(\cA,X)_H$. Since pull-backs and push-forwards are maps of $B(G)$-modules, they preserve this decomposition.\\

\subsection{Localization in modular K-theory}

Now, following the ideas expressed in the previous paragraph, we would like to compare the modular $K-$theory of $X$ with that of $X^H$ for any $H<G$. We do so by comparing the modular K-theory of $[X/G]$ with that of its wild inertia.

For any $\cA$ (finitely generated projectively complete full subcategory of $k[G]-\bf{mod}$ containing the trivial representation and equipped with an action of the Burnside algebra), we have an induction map $\alpha_\cA\colon K^*(\alpha^*\cA,X)\arr K^*(\cA,X)$ associated to the restriction $\alpha\colon k[G]-\bf{mod}\arr k[N_G(H)]-\bf{mod}$, and an induction $i_\cA\colon K^*(\alpha^*\cA,X^H)\arr K^*(\alpha^*\cA,X)$. Composing them we get an induction functor
\[
\hat i_\cA\colon K^*(\alpha^*\cA,X^H)\arr K^*(\cA,X).
\]
Our main result is the following localization theorem:

\begin{theorem}
Let $(X,\cA)$ be as above. Suppose that $X$ admits a stratification
\[
U_0\xrightarrow{j_0} U_1\xrightarrow{j_1}\cdots \xrightarrow{j_{n-1}} U_n=X
\]
by open $G-$invariant sets such that for any stratum $X_i:=U_{i+1}\backslash U_i$ the quotient stack $[X_i/G]$ is a gerbe with stabilizer normal in $G$.

Then for any $H<G$ the induction 
\[
\hat i_\cA\colon K^*(\alpha^*\cA,X^H)_H\arr K^*(\cA,X)_H
\]
is an isomorphism.
\end{theorem}

\begin{proof}
First of all we use that $X$ admits a stratification
\[
U_0\xrightarrow{j_0} U_1\xrightarrow{j_1}\cdots \xrightarrow{j_{n-1}} U_n=X
\]
by open $G-$invariant sets such that for any stratum $X_i:=U_{i+1}\backslash U_i$ the quotient stack $[X_i/G]$ is a gerbe with normal stabilizer. Using this fact, we can prove that the theorem holds for each $U_i$ by induction on $i$. Indeed, by Propositon~\ref{gerbe} any closed immersion $i_i\colon X_i\arr U_{i+1}$ satisfies the condition $(*)$ of Theorem~\ref{localization}, and so we have fon any $l$ a commuting diagram
\[
\begin{tikzcd}
K_{l+1}(\alpha^*\cA,U_i^H)\ar[r]\ar[d]&K_{l}(\alpha^*\cA,X_i^H)\ar[r,"i_{i,\cA}"]\ar[d]
&K_{l}(\alpha^*\cA,U_{i+1}^H)\ar[r,"j_i^{\cA}"]\ar[d] &K_{l}(\alpha^*\cA,U_i^H)\ar[r]\ar[d]& K_{l-1}(\alpha^*\cA,X_i^H)\ar[d]\\
K_{l+1}(\cA,U_i)\ar[r]&K_{l}(\cA,X_i)\ar[r,"i_{i,\cA}"]
&K_{l}(\cA,U_{i+1})\ar[r,"j_i^{\cA}"] &K_{l}(\cA,U_i)\ar[r]& K_{l-1}(\cA,X_i)
\end{tikzcd}
\]
From the five lemma we immediately infer that if we know the thesis for $U_i$ and for any $X$ such that $[X/G]$ is a gerbe then we also have it for $U_{i+1}$.\\

All we have to do, then, is to show the thesis in the case when $[X/G]$ is a gerbe. In that case there is a normal subgroup $Q<G$ such that $Q$ acts trivially on $X$ (so that $[X/G]$ is banded by $Q$) and $P:=G/Q$ acts freely on $X$. Let us consider the pull back $\alpha^*\cA$ with respect to $\alpha: k[G]-\bf{mod}\arr k[Q]-\bf{mod}$; note that $P$ has an action both on $X$ and on $\cA_{|Q}$, and so (by pull-back) it has an action on $\cA_X$ and on $K^*(\cA_{|Q},X)$.

Our first step will be to show that $\alpha^\cA$ induces an isomorphism
\[
\alpha^\cA\colon K^*(\cA,X)\xrightarrow{\simeq} K^*(\cA_{|Q},X)^P.
\]
To see it, we will show that both $\alpha_\cA\alpha^\cA$ and $\alpha^\cA\alpha_\cA$ are isomorphisms.

Consider the former. Let $Y:=X/P$; we have a commuting square
\[
\begin{tikzcd}
\lbrack X/Q \rbrack\ar[r]\ar[d] & \lbrack X/G \rbrack\ar[d,"\pi"]\\
X\ar[r,"\tilde\alpha"] & Y
\end{tikzcd}
\]
and we have $(\cA_{|Q})_X\simeq \tilde\alpha^*\cA_X$. In particular the composition $\alpha_\cA\alpha^\cA$ is just induced by the composition $\tilde\alpha_*\tilde\alpha^*$ on $Qcoh(Y)$. By the classical projection formula, this functor is just the tensor product with the locally free sheaf $\cE:=\tilde\alpha_*\cO_X=\pi_*s_X^*(Ind^G_Q1)$. By Morita theory this is a local progenerator, so it induces an isomorphism of the exact category $Coh(\cA,X)$ and thus an isomorphism in K-theory; alternatively we can observe that the tensor product with a locally free sheaf induces an action of $K^*(Y)$ on $K^*(\cA,X)$, and the class $\lbrack\cE\rbrack$ is a unit in $K_0(Y)$.

Now consider the latter composition. The push-forward $\alpha_\cA\cF$ of an $\cA_{|Q}-$sheaf $\cF$ sends $s_X^*V$ to the sheaf $\tilde\alpha_*\cF(V)$, with $(\cA_{|Q})_X$-action induced by the map $\cA_X\arr \tilde\alpha_*\tilde\alpha^*\cA_X$. Then composing with the pull-back we get 
\[
s_X^*V\arr\bigotimes_{p\in P} p^*F(V)^p 
\]
where $(-)^p$ means that the $\alpha^*\cA_X$-action is twisted by $p$. 
We conclude that $\alpha^\cA\alpha_\cA\cF=\sum_{p\in P}p\cdot\cF$; thus the functor $\alpha^\cA\alpha_\cA$ coincides with $\sum_{p\in P}p\cdot(-)$ in $Coh(\cA,X)$, and this relation must then hold also in K-theory. We conclude, since $\sum_{p\in P}p\cdot(-)$ is an isomorphism in $K^*(\alpha^*\cA,X)^P$.\\

Now, from this fact it follows trivially that we just need to consider the case when $H<Q$ (otherwise $K^*(\cA_{|Q},X)_H$ is zero and $X^H$ is empty); moreover we have $N_G(H)/N_Q(H)\simeq P$, so that also
\[
\alpha^\cA\colon K^*(\cA_{|N_G(H)},X^H)\xrightarrow{\simeq} K^*(\cA_{|N_Q(H)},X^H)^P
\] 
and we just need to prove the claim for $K^*(\cA_{|Q},X)$. In particular we can suppose that $G$ acts trivially and the Auslander algebra is constant on $X$.

In this case, the result follows from the following
\begin{lemma}\call{conlon}
Suppose that the Auslander algebra is trivial. Then the induction $\alpha_\cA$ and restriction $\alpha^\cA$ induce an isomorphism
\[
K^*(\cA,X)_H\simeq K^*(\cA_{|H},X)^{N_G(H)}.
\]
\end{lemma}

\begin{proof}
To avoid confusion, let us call $\hat i_\cA=\alpha_cA$ and $\hat i^\cA=\alpha^\cA$. Again, we would like to calculate the two compositions $\hat i_\cA\hat i^\cA$ and $\hat i^\cA\hat i_\cA$.

Let $I$ be the direct sum of a set of generators for $\cA$. Then the push-forward $\hat i_\cA$ corresponds to the push-forward from $End_H(I)-$modules to $End_G(I)$-modules, whereas the pull-back corresponds to $\otimes_{End(I)^G}End(I)^H$, where $End(I)^G$ acts on $\cF(I)$ on the left and on $End(I)^H$ on the right. The left $End(I)^H$-module structure on the pull-back is induced by the left action of $End(I)^H$ on itself.

Consider first $\hat i_\cA\hat i^\cA$, and let $T:=Ind^G_H(1)$, which by hypothesis is an element of $\cA$. By definition the sheaf $\hat i_\cA\hat i^\cA\cF$ sends an element $V\in \cA$ to the coend
\[
\int^W Hom_H(V,W)\otimes \cF(W).
\]
But $Hom_H(V,W)=Hom_H(Hom(V,W),1)=Hom_G(Hom(V,W),T)$, so $\hat i_\cA\hat i^\cA\cF(V)$ is equal to
\[
\int^W Hom_G(Hom(V,W),T)\otimes \cF(W)
\]
In particular, $\hat i_\cA\hat i^\cA\cF=\cF\otimes_{\cA_X}\cA_T$. We conclude that $\hat i_\cA\hat i^\cA$ is equal to the multiplication by $\cA_T$ in $K^*(\cA,X)_H$; this is an isomorphism, since $T$ is invertible in the $H-$part $B(G)_H$ (it is actually a multiple of the generator).

In order to deal with $\hat i^\cA\hat i_\cA$, we want to prove a Mackey-type formula. Clearly it is sufficient to prove it in the exact category $Coh(\cA_{H},X)$ and, following the definitions, for the representable sheaf $\cF=\cA_{I|H}$. Indeed, by hypothesis, we have $Res^G_HInd^G_H\cF=\cF\otimes _{End(I)^G}End(I)^H$, where $End(I)^H$ is seen as a right $End(I)^G-$module and as a left $End(I)^H$-module. But $\cF\otimes _{End(I)^G}End(I)^H=\cF\otimes _{End(I)^H}(End(I)^H\otimes _{End(I)^G}End(I)^H)$, where $End(I)^H\otimes _{End(I)^G}End(I)^H$ has a right $End(I)^H-$module structure given by the right action of itself on the first copy of $End(I)^H$.

Now, the induced representation $I':=Ind^G_H I_{|H}=I\otimes_G T$ is representable in $\cA$, so the push-forward sends $\alpha^*\cA_I$ to $\cA_{I'}$, and composing with the pull-back ge get $\alpha^*\cA_{Res^G_HInd^G_H I_{|H}}$. By the classical Mackey formula, this is equal to a sum $\underset{g\in H\backslash G\slash H} {\bigoplus}Ind^H_{^gH_i}Res^H_{H_i}(I_{|H})$. By Proposition~\ref{restriction} we can assume that each of the $I_i:=Ind^H_{H_i^g}Res^H_{H_i}(I_{|H})$ lies in $\cA_{|H}$; we see then that $$Res^G_HInd^G_H\cA_{I|H}(I_{|H})=\bigoplus Hom_H(I,I_i)=\bigoplus \alpha^*\cA_{I_i}(I_{|H}).$$
 By the same reasoning as above, we have that $$Ind^H_{H_i^g}Res^H_{H_i} \alpha^*\cA_{I|H}(I)= \alpha^*\cA_{I_i}(I),$$
 so the result holds, since this is an isomorphism of $End(I)^H-$bimodules. 

We saw that we have a Mackey-type formula
\[
Res^G_HInd^G_H\cF = \underset{g\in H\backslash G\slash H} {\bigoplus}Ind^H_{^gH_i}Res^H_{H_i} \cF
\]
in the category $Coh(\cA_{|H},X)$, and thus also in K-theory. However when $H_i=H\cap H^g$ is not equal to $H$, the restriction $Res^H_{H_i}$ is zero on the $H-$part, hence we have $Res^G_HInd^G_H\cF=\underset{n\in N_G(H)\slash H}{\bigoplus}n\cdot\cF$ (where, again, the $N_G(H)$-action is induced by the conjugation action on $\cA_{|H}$). However this functor is a multiple of the identity on $K^*(\cA_{|H},X)^{N_G(H)}$, and we conclude.
\end{proof}

Comparing the result of the Lemma applied to $G$ with the one applied to $N_G(H)$, (and noting that obviously $N_{N_G(H)}(H)=N_G(H)$) we immediately have our thesis.
\end{proof}

\begin{remark}
By following the proof of the Theorem, we can also prove that the induction
\[
\alpha_\cA\colon K^*(\cA_{|H},X^H)^{N_G(H)}\arr K^*(\cA_{N_G(H)},X^H)
\]
is an isomorphism.
\end{remark}

\begin{proof}
Exactly as in the proof of the preceding Theorem we use the existence of a stratification of the quotient $[X/G]$ by gerbes and the five lemma to reduce to the case where there exists $Q$ normal in $G$ acting trivially on $X$ and such that $P=G/Q$ acts freely on $X$.

But then we have a commuting diagram
\[
\begin{tikzcd}
K^*(\cA_{|H},X)^{N_G(H)}\ar[r,"\simeq"]\ar[d]&(K^*(\cA_{|H},X)^{N_Q(H)})^P\ar[d]\\
K^*(\cA_{|N_G(H)},X) \ar[r,"\simeq"] & K^*(\cA_{|N_Q(H)},X)^P
\end{tikzcd}
\]
so we can suppose that $H<Q$ or in other words that $G$ acts trivially on $X$.

In this case the result is exactly the content of Lemma~\ref{conlon} .
\end{proof}

\begin{remark}
Suppose that $G$ is of finite representation type and $\cA=k[G]-\bf{mod}$. Then $K^*(\cA-\bf{mod})\simeq K^*(\cA^{split})$.

By Theorem~\ref{semidirect}, in this case the only relevant subgroups are those of the form $P\rtimes C$, where $P$ is a $p-$group and $C$ is a cyclic group of order coprime with $p$ if $k$ is a characteristic $p$ field.
\end{remark}


\begin{thebibliography}{99}

\bibitem{bor} N. Borne, \textit{Cohomology of  $G$-sheaves in positive characteristic}, Adv. Math. 201 (2006), no. 2, 454-515.
\bibitem{conl} S. B. Conlon, \textit{Decompositions induced from the Burnside algebra}, J. Algebra 10 (1968), 102-122.
\bibitem{HKR} M. J. Hopkins, N. J. Kuhn, D. C. Ravenel, \textit{Generalized group characters and complex oriented cohomology theories}, J. Amer. Math. Soc. 13 (2000), no. 3, 553-594.
\bibitem{May} L. G. Lewis, J. P. May, M Steinberger, \textit{Equivariant Stable Homotopy Theory}, Springer, Lecture Notes in Mathematics, volume 1213.
\bibitem{Nak} S. Nakajima, \textit{Action of an automorphism of order  $p$  on cohomology groups of an algebraic curve}, J. Pure Appl. Algebra 42 (1986), no. 1, 85-94.
\bibitem{To} B. Toen, \textit{Théorèmes de Riemann-Roch pour les champs de Deligne-Mumford}, $K $-Theory 18 (1999), no. 1, 33-76.
\bibitem{Vistoli} A. Vistoli, \textit{Higher equivariant  K -theory for finite group actions}, Duke Math. J. 63 (1991), no. 2, 399-419.

\bibitem{yas} T. Yasuda, \textit{Motivic integration over wild Deligne-Mumford stacks}, Algebr. Geom. 11 (2024), no. 2, 178-255.

\end{thebibliography}
\end{document}